# A mathematical model for understanding and controlling monkeypox transmission dynamics in the United States and its implications for future epidemic management

Md. Azmir Ibne Islam[a,1], M H M Mubassir[a,b,2], Arindam Kumar Paul[c,3], Sharmin Sultana Shanta[c,4]

[a]*Department of Mathematics and Natural Sciences, Brac University, Dhaka 1212, Bangladesh*

[b]*Institute of Bioinformatics, University of Georgia, USA*

[c]*Mathematics Discipline, Khulna University, Khulna 9208, Bangladesh*

ARTICLE INFO

ABSTRACT

**Background:** Although the outbreak of human monkeypox (mpox) caused by the monkeypox virus (MPXV) has slowed down around the world, little is known about the short-term dynamics of this disease. This limited information highlights the critical need to assess the underlying interventions.

**Method:** To identify and re-examine the key pattern of the disease, a modified logistic growth model is presented and analysed in this paper. Our main focus is on the two non-pharmaceutical interventions: policies aimed at reducing human-to-human transmission and animal-to-human transmission. We incorporated these two strategies in the model as control parameters to understand their short-term significance on the epidemic, and to analyse their strengths in minimizing the infected cases. We used mpox data set of the United States from 10 May 2022 to 31 December 2022 in the model and estimated the baseline parameters.

**Results:** The model reveals a complying acceptance to the US data set. Model simulations highlight that preventive measures could play important roles in controlling the deadly spread of the disease in the year of 2022. During the transmission period, better outcomes could have been possible to achieve in the US if both controls were brought to action simultaneously.

**Conclusion:** Our study reflects that continuous application of the preventive strategies might be an effective tool to prevent the short-term outbreak of mpox or similar diseases. Moreover, such strategies could play supporting roles during pre- and post-vaccination periods.

[1] Corresponding author: Email: azmir.islam@bracu.ac.bd (MAII); ORCID: 0000-0002-7076-1636 (MAII)
[2] Email: mm07305@uga.edu (MHMM); ORCID: 0000-0002-9674-2425 (MHMM)
[3] Email: arindam017@gmail.com (AKP); ORCID: 0000-0003-3540-6919 (AKP)
[4] Email: shanta@math.ku.ac.bd (SSS); ORCID: 0000-0002-8289-6749 (SSS)

## 1. Introduction

For over two years, the world has suffered from the deadly viral disease COVID-19 (SARS-CoV-2). Before this pandemic ended, another old zoonotic viral disease known as monkeypox (mpox) came to focus, which caused significant cases in 2022 [14, 24, 39]. The mpox disease is caused by the monkeypox virus (MPXV), an Orthopoxvirus genus virus. It was first discovered in cynomolgus monkeys in 1958 during a study of the polio vaccine in Copenhagen, Denmark [21, 27]. The first human mpox disease was reported in August 1970 in Bokenda, Democratic Republic of the Congo (DRC) [5, 12]. Then, through time, it became endemic in the African region [12, 26, 31].

The first mpox outbreak outside of Africa was reported in the US in 2003. According to the Centers for Disease Control and Prevention (CDC), 72 human mpox probable cases were registered in 2003 in the US due to close contact with pets [6, 17, 30]. The UK Health Security Agency (UKHSA) documented the first non-endemic mpox case of 2022 on May 7, 2022, and the number of cases had increased rapidly after that period [36]. On June 18, the Massachusetts Department of Public Health confirmed the first case of MPXV infection in the state and the first confirmed mpox case in the United States in 2022 [23]. Following then, the number of mpox cases in the US climbed dramatically. However, the current trend of mpox in the US is seen to be at a stable position (low level), but nearly around 30 thousand verified cases have been documented within 31 December 2022 due to this disease. As of 1 September 2023, more than 89 thousand confirmed cases have been reported around the world due to mpox [6].

Although MPXV is typically transmitted through animal bites or interactions, cases of 2022 outbreaks reveal varied transmission mechanisms [2, 37]. In 2022, most cases were discovered in non-endemic locations, but none of the patients had been to Africa. Human-to-human transmission plays a crucial role in the transmission of MPXV in such instances, which raises concern following the SARS-CoV-2 pandemic. There is also speculation that mpox is a sexually transmitted virus that can spread among guys who typically have sex with men (MSM) [13, 33]. Besides human-to-human transmission, animal-to-human transmission is also responsible for the outbreak of mpox [6]. The major signs of the disease include fever, rash and lymphadenopathy, which can lead to death if it becomes complex [19, 22, 25]. Other symptoms, including rectal pain and bleeding, were noted in the new outbreak for the US patients that had not previously been observed [10]. JYNNEOS is the first mpox vaccine licensed by the Food and Drug Administration (FDA) but is only available to people at risk of exposure [32].

Epidemiologists widely employ mathematical models as crucial tools to gain understanding of the spatial and temporal patterns of infectious diseases [3]. Models are used to assist health authorities and governments in deciding what preventive measures are needed to put in place and how to properly allocate the health facilities with limited resources [1]. Mathematical modelling is the technique to understand and/or forecast the dynamics of a disease [4, 11, 18].

The current work is devoted to the well-known logistic growth model which was first introduced by P. F. Verhulst (1804-1849), a Belgian mathematician and biologist, who was interested in forecasting the human population of numerous nations around 1840 [38]. In this paper, a modified logistic growth model with two interventions is considered so that the short-term dynamics of mpox transmission in the United States (during 2022) can be understood and analysed. The paper aims to re-examine the outbreak of mpox in the US from 10 May 2022 to 31 December 2022 based on the daily and the cumulative cases data. The ultimate focus of this work is on the effectiveness of the incorporated interventions and their implications towards future epidemics. The remainder of the current paper is organised as follows: section 2 introduces the model, shows properties of the controls, and considers data fitting. Section 3 discusses numerical results. Finally, section 4 concludes with limitations.

## 2. Methods

### *2.1. Model description*

The logistic model, first developed by P. F. Verhulst, has become a key tool in biology, epidemiology, economics, and finance [9]. Specifically, it has been applied to model various phenomena, such as population growth and bacterial proliferation in a culture [34]. The well-known logistic equation is defined as

$$\dot{y} = ry\left(1-\frac{y}{k}\right) \tag{2.1}$$

where $r$ is the growth rate, $k$ is the carrying capacity, and $y$ or $y(t)$ is the solution.

In recent years, researchers have used model (2.1) in understanding the dynamics of COVID-19 [7, 29, 35]. Since logistic equation (2.1) is frequently used to model population growth, it can be applied in modelling mpox dynamics too. Different types of modifications in the logistic model have been seen in recent works [15, 16, 20]. With this view, we propose the following modified model to study the short-term dynamics of mpox transmission:

$$\dot{y}(t) = (1-u)r\,y(t)\left(1-\frac{y(t)}{1/\gamma}\right)+(1-v)h \qquad (2.2)$$

In model (2.2), $y(t)$ denotes the cumulative number of infected cases due to the MPXV at time $t$. Parameters $r$ and $h$ denote human-to-human transmission rate, and zoonotic (animal-to-human) transmission rate, respectively. The carrying capacity $k$ is replaced by $1/\gamma$, which represents the final size of the epidemic and the parameter $\gamma$ represents the treatment facility rate. The parameters of model (2.2) are all positive (biologically meaningful), and listed in Table 1. Model (2.2) has two controls $u$ and $v$ (see subsection 2.2 for details). One assumption of this model is: in absence of any non-pharmaceutical interventions (meaning, $u$, $v = 0$) during the disease progression, the treatment facility (i.e., $\gamma$) provided by the government by default acts as a control policy. Thus, in the baseline scenario, we set $u$ and $v$ to be zero.

Table 1: Parameters of model (2.2).

| **Parameter** | **Definition** |
|---|---|
| $r$ | human-to-human transmission rate |
| $\gamma$ | treatment facility rate |
| $h$ | zoonotic (animal-to-human) transmission rate |

*2.2. Short-term control-design*

Two inputs $0 \le u < 1$ and $0 \le v < 1$ are introduced in model (2.2) as control parameters, which have influence on the transmission rates $r$ and $h$, respectively. Control 1 (i.e., $u$), incorporated in model (2.2), mainly signifies the strength of the interventions taken by the government for minimizing the human-to-human transmission rate, whereas control 2 (i.e., $v$) reflects the strength of the current interventions for minimizing the animal-to-human transmission rate. The values of $u$ and $v$ closer to 0, indicate the lower efficacy of control measures, whereas values closer to 1, signify higher efficacy. We assume the upper bounds of $u$ and $v$ to be less than 1, because it is hardly possible to receive complete efficacy levels due to certain constraints. To clearly understand the effectiveness of the controls in mitigating the spread of monkey pox virus, following strategies are considered: (i) Strategy 1: Implementation of control 1 only, (ii) Strategy 2: Implementation of control 2 only, and (iii) Strategy 3: Implementation of both controls simultaneously.

*2.3. Model analysis*

Model (2.2) is straightforward, which has the following exact solution

$$y(t)=\frac{1}{2}\left[A_2-\frac{B_1}{B_2}\tanh\left(\frac{1}{2}\left(-B_1B_3c_1-B_1B_2t\right)\right)\right] \tag{2.3}$$

where the constants are

$$c_1, A_1=(1-u)r, A_2=\frac{1}{\gamma}, A_3=(1-v)h,$$

$$B_1=\sqrt{A_1A_2+4A_3},\ B_2=\sqrt{\frac{A_1}{A_2}},\ B_3=\sqrt{A_1A_2}.$$

Note that, $c_1$ is the constant of integration. Appendix A provides more details regarding this exact solution.

We also compute the equilibrium points (fixed points) of model (2.2) by considering $\dot{y}=0$. The equilibria of model (2.2) are

$$\xi_1=\frac{A_1A_2+A_2B_1B_2}{2A_1} \text{ and } \xi_2=\frac{A_1A_2-B_1B_3}{2A_1}.$$

We verified numerically that equilibrium point $\xi_2$ gives a negative result for the baseline parameters (computed by data fitting in subsection 2.4). Since $\xi_2$ corresponds to biologically non-viable scenarios, we ignored it. We thus only consider the equilibrium point $\xi_1$ in our analysis, which is positive and physically meaningful.

*2.4. Data*

Model (2.2) is fitted to the US mpox data set [6] in order to estimate the three parameters: $r$, $\gamma$ and $h$. Starting with some initial guesses for these parameters, we solved model (2.2) using the built-in 'ode45' function of MATLAB. We then applied the resulting solutions to MATLAB routine 'fmincon' to estimate these parameters that minimize the error function $E$; for data fitting, the error function $E$ is minimized as follows

$$E=\sum_{i=1}^{236}\left(y(t_i)-\hat{y}(t_i)\right)^2 \tag{2.4}$$

where $y(t_i)$ represents simulation results and $\hat{y}(t_i)$ denotes the data of US cumulative cases. Moreover, $i$ = 1 indicates 10 May 2022, whereas $i$ = 236 signifies 31 December 2022.

We also considered Monte Carlo simulation to account the uncertainty in the estimates of $r$, $\gamma$ and $h$. We used the estimates found from 'fmincon' as the baseline and introduced random sampling across plausible parameter ranges. We assumed that each parameter ($r$, $\gamma$ and $h$) follows a uniform distribution in the predefined range. We then performed 10000 Monte Carlo simulations in order to estimate the distribution of each parameter. Through this approach, we captured the variability and uncertainty in the parameter estimates. The final values of the parameters along with their corresponding 95% confidence intervals are provided in Table 2; we also refer Figure 1 for details.

Table 2: Fitted values of the parameters of model (2.2).

| **Parameter** | **Fitted value** | **Range** (95% confidence interval) | **Unit** |
|---|---|---|---|
| $r$ | 0.06 | [0.0509, 0.0651] | per day |
| $\gamma$ | 0.000034 | [0.000032, 0.000036] | per day |
| $h$ | 5.99 | [5.9712, 6.0134] | per day |

*2.5. Sensitivity analysis*

Following [8, 28], we performed sensitivity analysis of the equilibrium point $\xi_1$ to assess how it changes in response to variations in the parameters. The sensitivities of $\xi_1$ with respect to $r$, $\gamma$ and $h$ are $\frac{\partial \xi_1}{\partial r}\frac{r}{\xi_1}$, $\frac{\partial \xi_1}{\partial \gamma}\frac{\gamma}{\xi_1}$ and $\frac{\partial \xi_1}{\partial h}\frac{h}{\xi_1}$, respectively. Table 3 represents the sensitivity indices of $\xi_1$ for the fitted values of the parameters given in Table 2. We found that treatment facility rate $\gamma$ is the most influential parameter because it has a strong negative elasticity, which indicates that small changes of $\gamma$ will have a significant impact on reducing $\xi_1$ (see Figure 6). The parameter $r$ has also a negative sensitivity index, but it is minimal compared to $\gamma$. On the contrary, the sensitivity index of $h$ is positive, meaning that, $\xi_1$ will slightly increase if $h$ increases.

Table 3: Sensitivity of $\xi_1$.

| **Parameter** | **Sensitivity index** |
|---|---|
| $r$ | -0.00336 |
| $\gamma$ | -0.99664 |
| $h$ | +0.00336 |

## 3. Results and Discussion

### *3.1. Mpox in the US*

Model (2.2) fits well with the US mpox daily and cumulative confirmed cases data from 10 May 2022 to 31 December 2022 (see Figure 1). It is seen that the peak of infection occurred in the middle of August 2022 with around 15 thousand cumulative cases (approx.). The simulation of model (2.2) depicts that there were nearly 450 patients at the epidemic peak. By the end of 2022, around 30 thousand population became infected by MPXV.

### *3.2. Effect of the implementation of control policy 1*

If it was possible to apply control policy 1 with 40% effectiveness, the cumulative cases could be minimized by 38% on average (see Figure 2). Due to this, there could be a delay in peak time (end of September 2022) with a peak value smaller than the baseline case. The final size could be brought down to more than 78% within the end of 2022 if control policy 1 was implemented with an efficacy level of 80%.

### *3.3. Effect of the implementation of control policy 2*

From Figure 3, it is seen that an average of 17% reduction of cumulative cases could be possible compared to the baseline case if the effectiveness of control policy 2 was 40% during the epidemic, whereas the total number of cases could be averted up to 40% if control policy 2 was 80%. Despite the variation in the efficacy level of the applied control, no significant change was seen other than the delay in peak timing.

### *3.4. Combined effect of the controls*

Due to the implementation of both controls at the same time, significant reductions in daily cases and cumulative cases are observed (see Figure 4). If both controls were applied at 40% effectiveness, around 52% decrease in cumulative cases could be seen. On an average of 95% decrease in the final size could be possible if both controls were implemented with 80% efficacy level. While applying both controls at different efficacy levels, delays in the peak of infection might be observed with less patients at the peak time compared to the baseline case.

### *3.5. Comparison between the controls*

Comparison between the two intervention policies is shown in Figure 5. In case of single control policy, control policy 1 could be better than control policy 2 in minimizing the average number of cumulative cases from the baseline case. The total number of patients could significantly be reduced at a greater scale if both controls were implemented at the same time. Therefore, the size of the epidemic peak could possibly be shorten down with less number of patients than that of observed from the US mpox data set [6].

### 3.6. Effect of treatment facility rate

Figure 6 represents that with the enhancement of treatment facilities, the final size of the epidemic could be reduced. If the treatment facilities were enhanced by two times, the cumulative cases could be minimized up to 31% compared to the base scenario. Epidemic size could also be curtailed down to 55% with the increase of treatment facility by five times. Due to the augmentation in the treatment rate, the peak time might occur earlier with significant smaller peaks compared to the current peak.

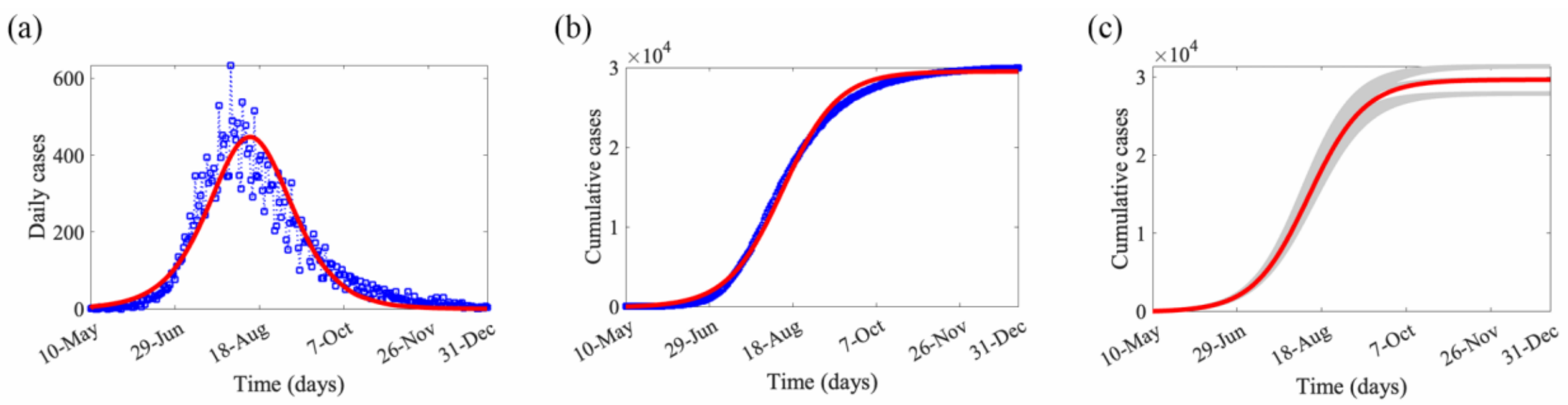

**Figure 1:** Mpox in the US from 10 May 2022 to 31 December 2022 (baseline case). Model (2.2) was fitted to the US mpox data [6]. Panel (a) represents model fitting (red) with US daily cases data (squares with dotted line), whereas panel (b) is the model fitting (red) with US cumulative cases data (squares with dotted line). Panel (c) shows model fitting (red) with uncertainty bounds (gray region) with 95% confidence intervals for the parameters; the red curves in panels (b) and (c) represent the mean trajectories. In panel (a), we calculated the daily cases from the cumulative cases data by computing the difference between successive cumulative values.

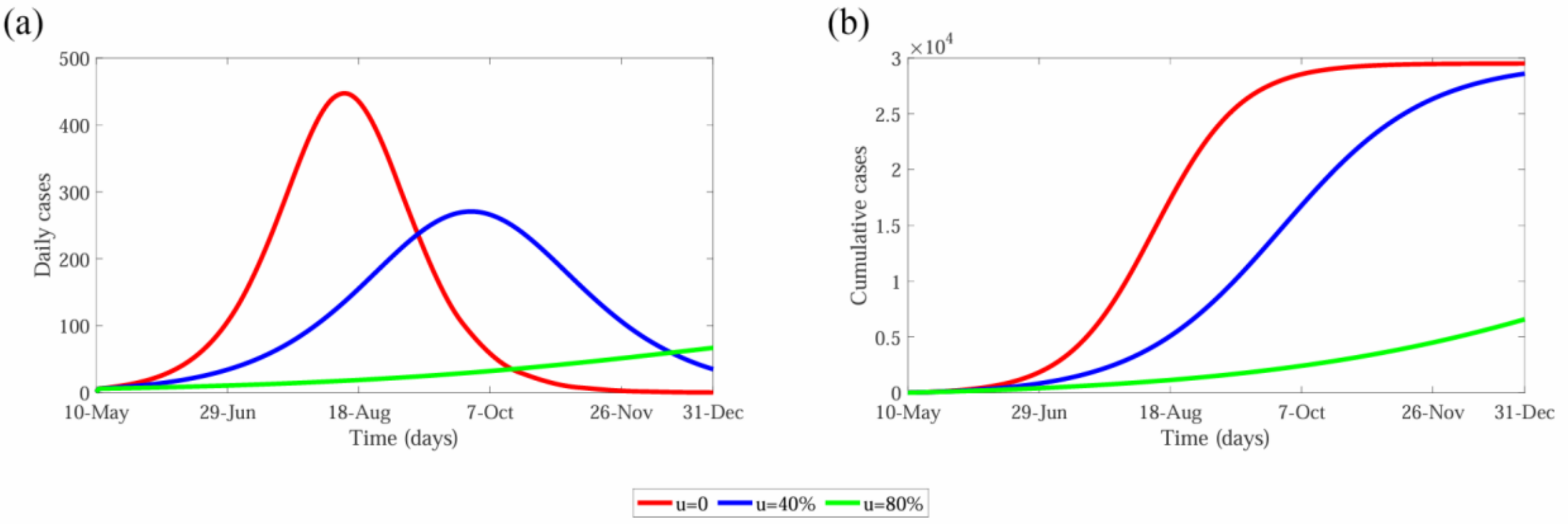

**Figure 2:** Effect of control policy 1 for a short-term period (from 10 May 2022 to 31 December 2022) for $u = 0$ (baseline in red), 0.4 (in blue) and 0.8 (in green). Panel (a) shows curves of daily cases, whereas panel (b) shows curves of cumulative cases. Fitted values of $r$, $\gamma$ and $h$ from Table 2 are used here with $v = 0$; (see Figure B.1(a), which shows full interactive effects of this analysis).

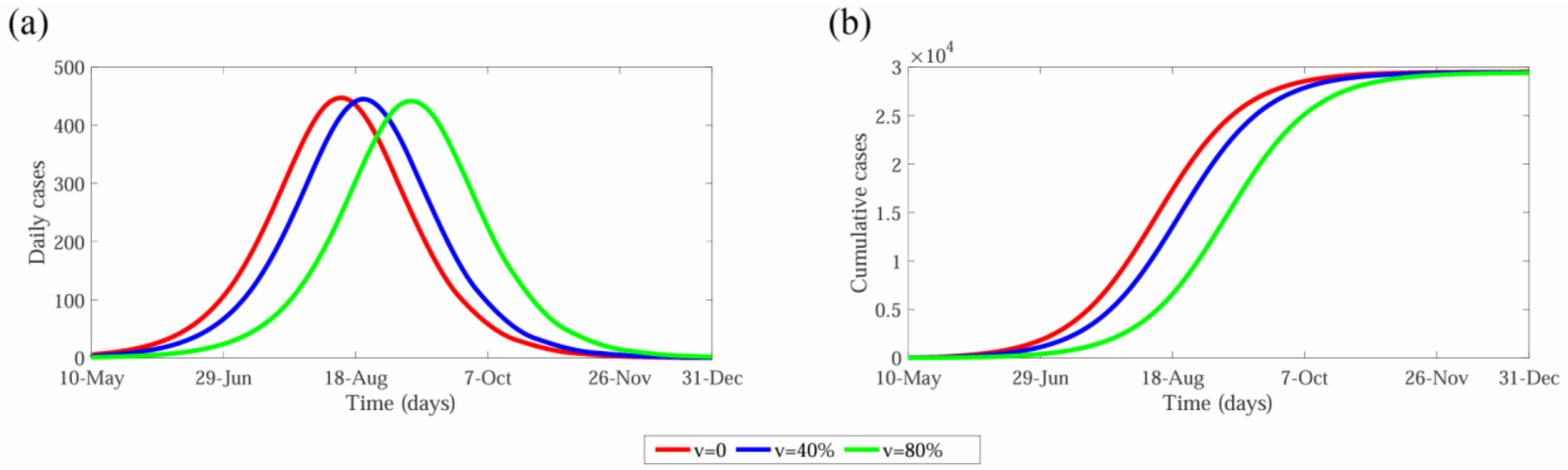

**Figure 3:** Effect of control policy 2 for a short-term period (from 10 May 2022 to 31 December 2022) for $v = 0$ (baseline in red), 0.4 (in blue) and 0.8 (in green). Panel (a) shows curves of daily cases, whereas panel (b) shows curves of cumulative cases. Fitted values of $r$, $\gamma$ and $h$ from Table 2 are used here with $u = 0$; (see Figure B.1(b), which shows full interactive effects of this analysis).

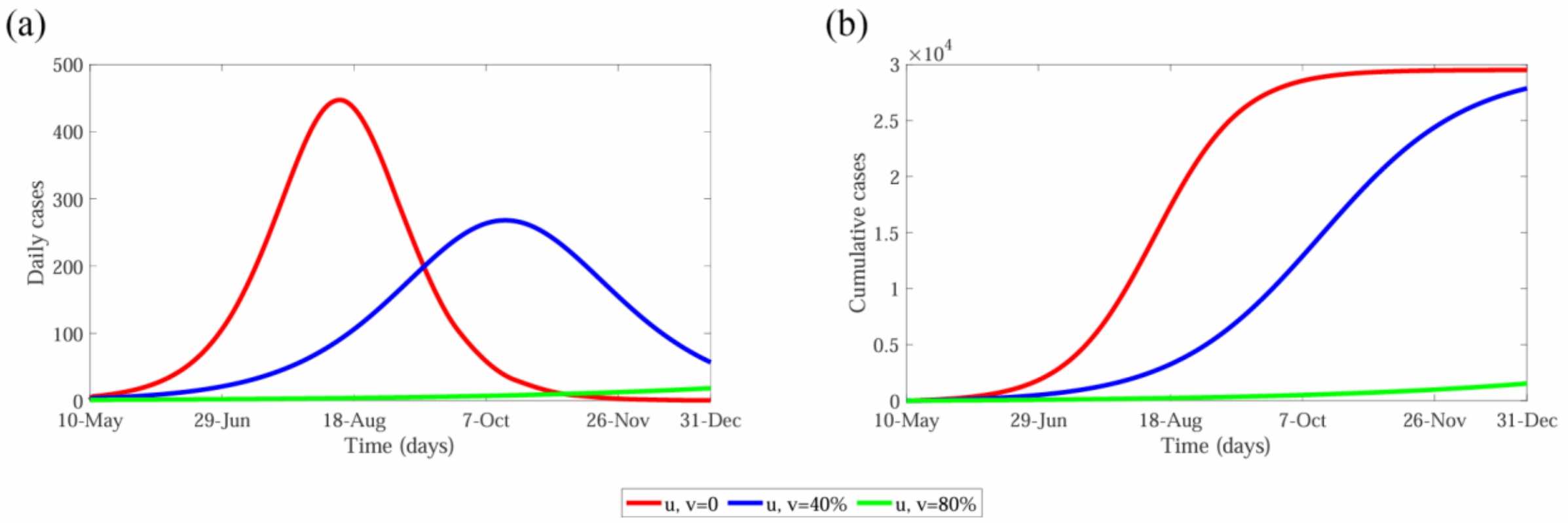

**Figure 4:** Combined effect of control policies 1 and 2 for a short-term period (from 10 May 2022 to 31 December 2022) for $u$, $v = 0$ (baseline in red), 0.4 (in blue) and 0.8 (in green). Panel (a) shows curves of daily cases, whereas panel (b) shows curves of cumulative cases. Fitted values of $r$, $\gamma$ and $h$ from Table 2 are used here; (see Figure B.1(c), which shows full interactive effects of this analysis).

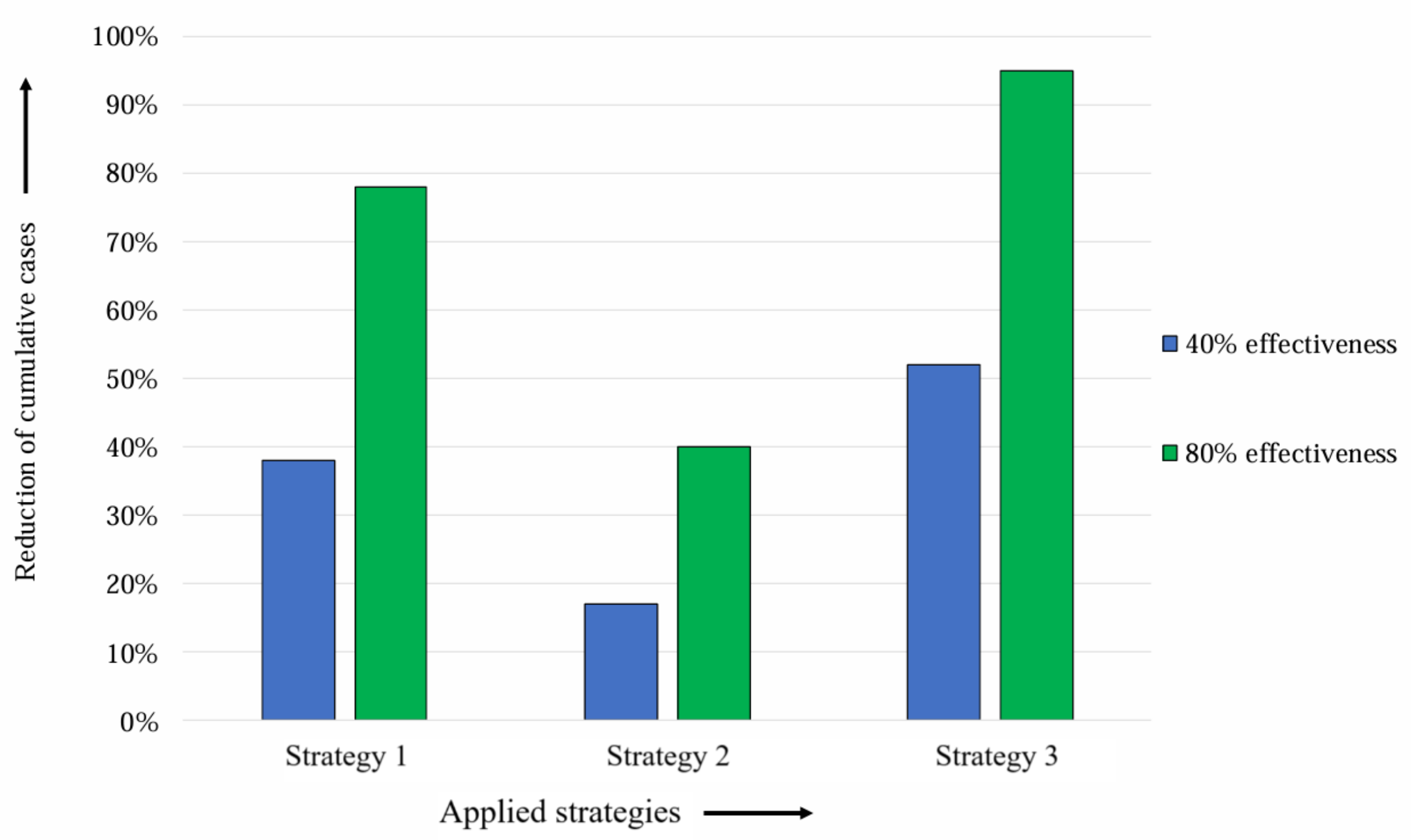

**Figure 5:** Comparison of the applied controls or strategies. The horizontal axis represents the strategies, and the vertical axis represents the percentage of reduction of cumulative cases. The two bars in the left reveal Strategy 1; the two bars in the middle are for Strategy 2; and the two bars in the right show Strategy 3.

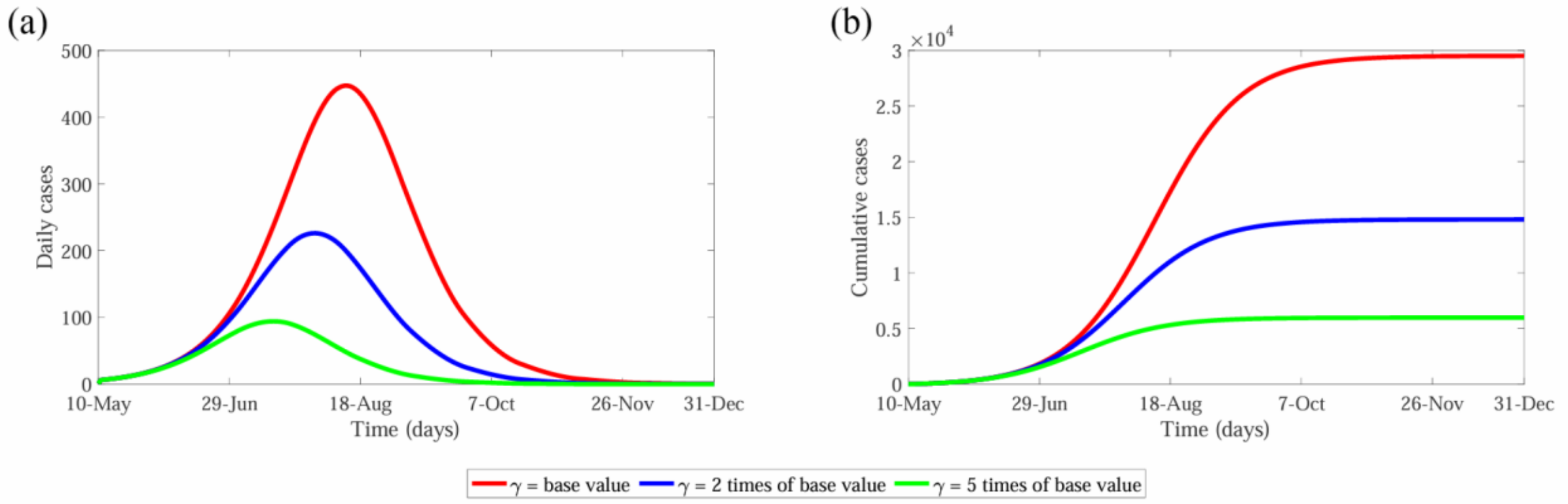

**Figure 6:** Effect of treatment facility rate (from 10 May 2022 to 31 December 2022) for $\gamma = 0.000034$ (baseline in red), $\gamma = 0.000068$ (in blue) and $\gamma = 0.00017$ (in green), with $u$, $v = 0$. Panel (a) shows curves of daily cases, whereas panel (b) shows curves of cumulative cases. Fitted values of $r$ and $h$ from Table 2 are used here; (see Figure B.1(d), which shows full interactive effects of this analysis).

## 4. Conclusion

In this paper, we present a modified logistic model with control strategies. Our model captures the short-term dynamics of mpox infection in the United States from 10 May 2022 to 31 December 2022. While re-examining the mpox trend in this observation period, it is revealed that control policies with a sound efficacy level could weaken the transmission flow and curb the final epidemic size. The timely implementation of preventive measures could have benefits in waning the peak size. The current study also reveals that reduction of human-to-human transmission should receive much more importance than minimizing animal-to-human transmission. Therefore, focus needs to be given more on the interventions that are responsible for the decline of human-to-human transmission.

Although our model shows a good agreement with the short-term data set, it has limitations. The results found from the model may possess a trivial variation from the real scenario. With large data set, more realistic outcomes could be unveiled. Moreover, the present model works with two constant control policies. However, to better understand the dynamics of the disease and for reasonable results, the model can be modified with optimal control strategies. Lastly, it is undeniable that our model is not based on long-term dynamics, and it might tend to infinity if parameters are relaxed. Considering such complexities in mind, the present model could be reconstructed to a nonautonomous dynamic model, not in terms of controls directly but in terms of the role of some external forces to restrict the transmission. Although a complete explanation of this concept is beyond the scope of the present paper, it sparks numerous questions and sets the platform to the readers for future exploration and research. However, from our side, we present *parameter continuation* of the model in Appendix C.

---

**Conflicts of interest:** The authors declare that there have no conflicts of interest.

**Ethical approval:** The study does not require any ethical approval.

**Use of AI tools:** The authors declare that no AI is used in the text.

**Data sharing statement:** All the data presented in the manuscript are publicly available.

**Funding:** There was no funding associated with this study.

**Author contribution:** MAII and MHMM designed and planned the study. MAII, AKP and SSS did the mathematical analysis. MHMM reviewed the manuscript. All authors approved the final version of the manuscript for publication.

**Supplementary:** https://github.com/azmiribneislam/Mpox_2024_MAII_MHMM_AKP_SSS

**Acknowledgements:** We, the authors, are grateful to Md. Saalim Shadmaan for thoroughly reading the paper multiple times while preparing this work; to Mohamed Bakheet for tremendous feedback and for reminding us about important mathematical analysis during the very first draft; to Raghib Ishraq Alvy, Muhammad Mahfuz Hasan, Arif Hossain Mazumder, Anika Ferdous, Md. Sajeebul Islam Sk., and Md. Rafsan Islam Rim for their insightful conversations; and to Md. Mahfuzur Rahman and Moumita Datta Gupta for sharing their statistical knowledge during the second draft. We also extend our gratitude to our parents for their unwavering support; we are indebted to our family members and children for giving us space and time to work. We especially thank Ashrafur Rahman for all the insightful online sessions on dynamical systems during the COVID-19 period. Finally, we thank the anonymous reviewers and the editor for their careful reading of all drafts and for the constructive suggestions.

**Dedication:** While this paper was under review, we tragically lost our co-author, Arindam Kumar Paul (1995 - 2024), in a road accident on April 27, 2024. We dedicate this paper to his memory. His invaluable dedication and contributions to this work will always be remembered.

---

## Appendix

### A. Exact solution of model (2.2)

In subsection 2.3, we present the exact solution of model (2.2). In this appendix, we determine the value of $c_1$ (constant of integration) and compare this exact solution with the numerical solution (obtained for the baseline case; see Figure 1(b)).

At $t = 0$ (corresponding to 10 May 2022), we have $y(0) = 1$ [6]. Using this initial condition along with the fitted values of $r$, $\gamma$ and $h$ from Table 2, we get $c_1 = -0.0032$. Thus, equation (2.3) gives a particular solution. Figure A.1(a) shows the comparison between the numerical curve and the exact curve. Both curves completely overlap on top of each other. This completely verifies our numerical findings. Figure A.1(b) shows a 3-dimensional plot of the exact solution for some choices of $c_1$, ranging from $-0.0050$ to $0.0050$. The family of solutions in $(c_1, t, y)$-space reveals how $y(t)$ behaves with the change of $c_1$ values. See the supplementary link for visualisation: https://github.com/azmiribneislam/Mpox_2024_MAII_MHMM_AKP_SSS.

### B. Short-term dynamics of model (2.2) for multiple values of the parameters

This appendix provides an additional visualisation of the results presented in section 3. Here, we present simulations of model (2.2) for a wide range of parameter values (for $u$, $v$ and $\gamma$). Effects are shown from 10 May 2022 to 31 December 2022. Figure B.1 shows the complete short-term dynamics of model (2.2) in an interactive way. Consistency in the results is observed for multiple values of the parameters. For example, cumulative case decreases with the increase of the parameters that we are interested in. Animations of the panels in Figure B.1 can be found in the supplementary link.

### C. Family of equilibria

In this appendix, we relax all the parameters ($r$, $\gamma$ $h$, $u$ and $v$) of model (2.2), which means, there is no restriction in taking the parameter values. We intend to capture the global picture of model (2.2) from a *dynamical systems* point of view; basically, we are interested in observing how the model behaves with the change of any parameters. We start with equilibrium point $\xi_1$ and compute all one-dimensional family of equilibria of the model for each of the individual parameters; note that, at baseline case, $\xi_1 = 29511.26$. Figure C.1(a-e) shows all numerical continuations from this value of $\xi_1$, which clearly reveals the dependency of model (2.2) on the parameters. An important aspect is seen in Figure C.1(e), which shows continuation of family of equilibria of the model for the parameter $v$. There is a saddle-node bifurcation (fold) at $v = 74.65$. At this value of $v$, the stable branch of the fixed point $\xi_1$ collide with the unstable branch of the other fixed point $\xi_2$. In subsection 2.3, we ignored $\xi_2$ because it is negative in the baseline case (i.e., $\xi_2 = -99.5$). However, from Figure C.1(e), it's now clear that $\xi_2$ can also be positive with the increase of $v$, and the only attractor is the branch of $\xi_1$. Note that, after $v = 74.65$, model (2.2) has no equilibria since the two equilibria of the system ($\xi_1$ and $\xi_2$) collide and annihilate each other.

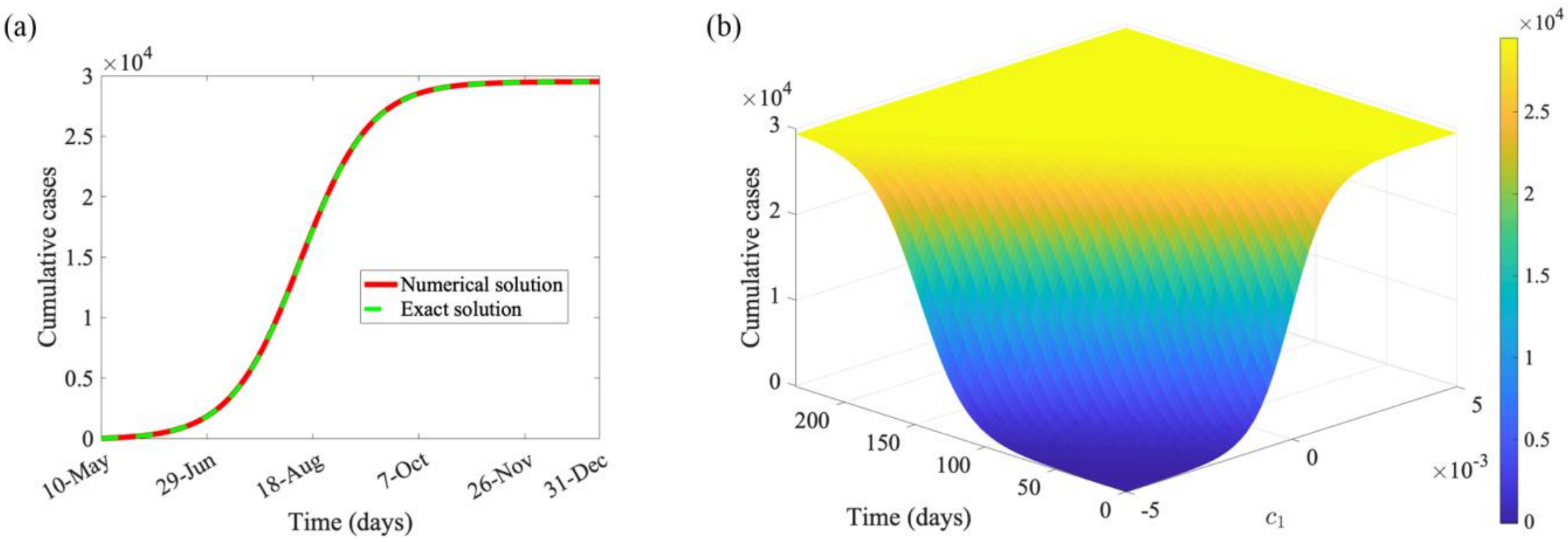

**Figure A.1:** Panel (a): exact solution versus numerical solution of model (2.2). Fitted values of $r$, $\gamma$ and $h$ from Table 2 are used here, with $u$, $v = 0$, and $c_1 = -0.0032$. Exact solution is shown in dashed line (green), whereas numerical solution is in solid line (red). Panel (b): family of solutions (exact) in $(c_1, t, y)$-space. Fitted values of $r$, $\gamma$ and $h$ from Table 2 are used here, with $u$, $v = 0$. The colour bar shows the values of $y$ (cumulative cases). In both panels, time frame is from 10 May 2022 to 31 December 2022; (see supplementary link for visualisation).

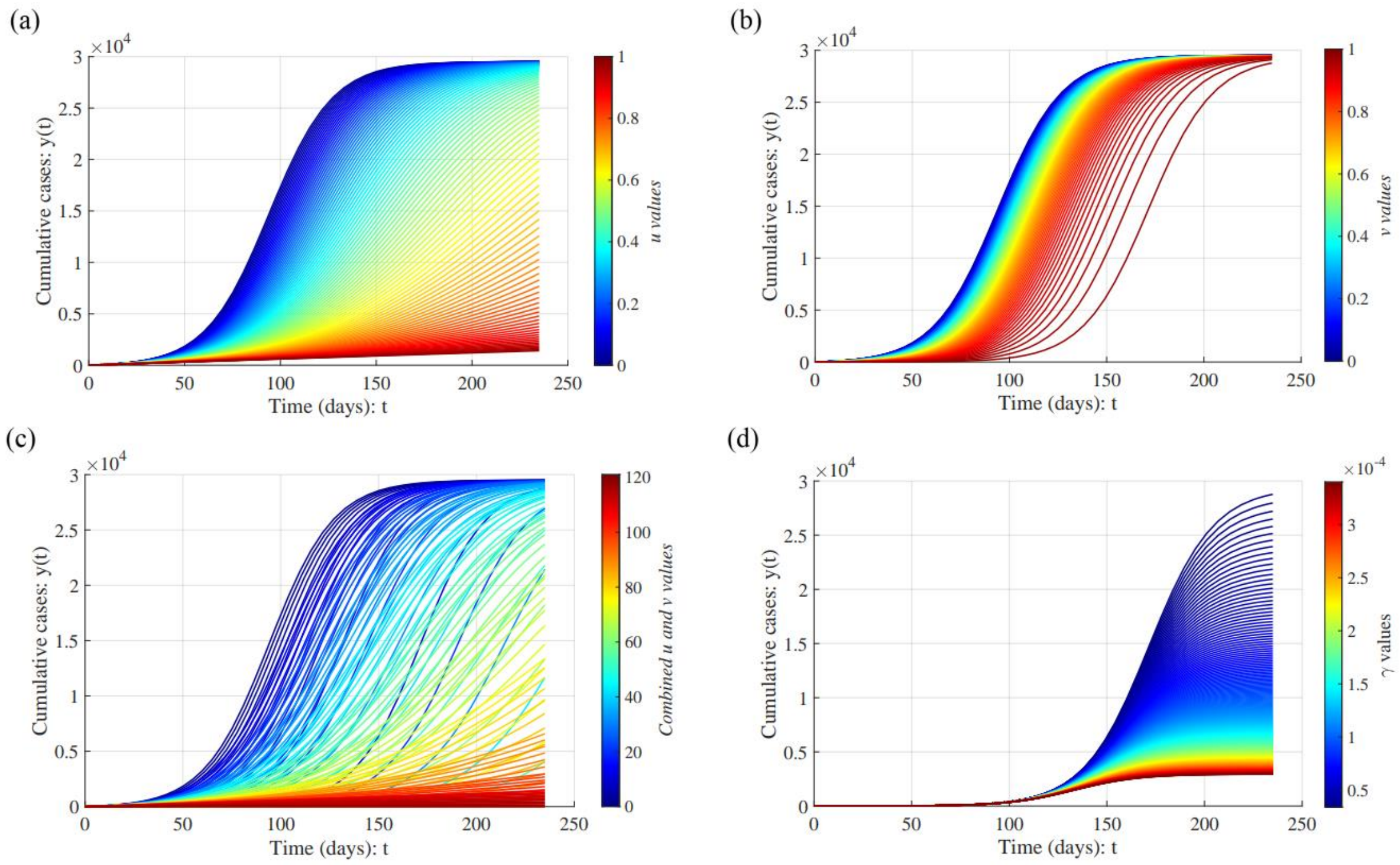

**Figure B.1:** Full interactive effects of: control policy 1 in panel (a); control policy 2 in panel (b); both control policies in panel (c); and treatment facility rate in panel (d). The colour bar in each panel shows the values of the respective parameter(s). In all panels, time frame is from 10 May 2022 to 31 December 2022; (see supplementary link for visualisation).

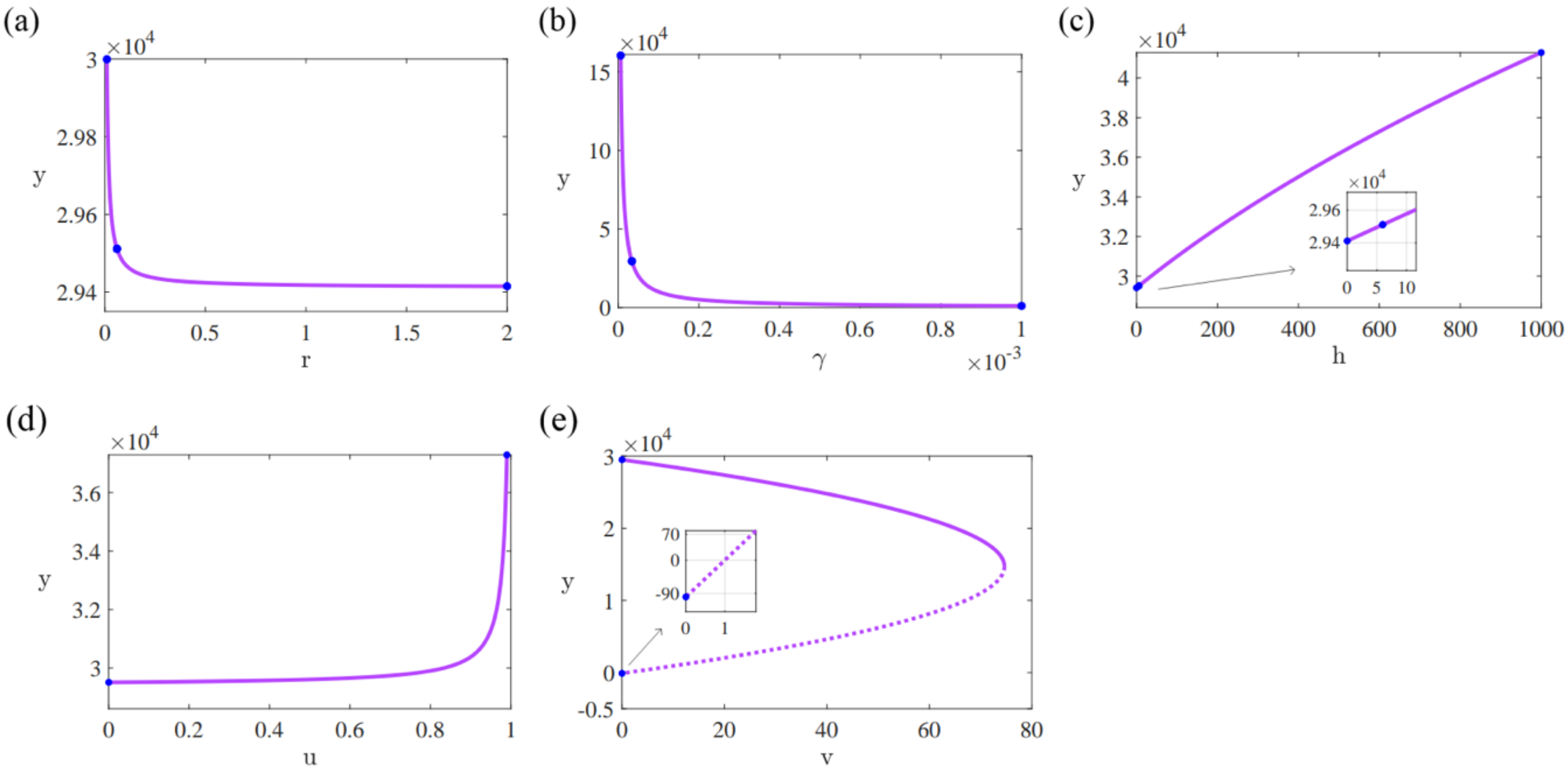

**Figure C.1:** Parameter continuation of model (2.2), showing family of equilibria: as a function of $r$ in panel (a); as a function of $\gamma$ in panel (b); as a function of $h$ in panel (c); as a function of $u$ in panel (d); and as a function of $v$ in panel (e). End points of each computation are marked via blue dots. Panels (c) and (e) show enlargements of tiny regions. In all panels, the solid purple curves are stable branches; but in panel (e), besides the stable branch, there is an unstable branch (dotted line in purple).